\documentclass[10pt]{article}
\usepackage{amsfonts,amsmath,amssymb,amsthm,cite}
\newtheorem{lem}{Lemma}
\begin{document}
\title{A note on parameter derivatives of classical orthogonal 
polynomials}
\author{Rados{\l}aw Szmytkowski \\*[3ex]
Atomic Physics Division, 
Department of Atomic Physics and Luminescence, \\
Faculty of Applied Physics and Mathematics, 
Gda{\'n}sk University of Technology, \\
Narutowicza 11/12, PL 80--233 Gda{\'n}sk, Poland \\
email: radek@mif.pg.gda.pl}
\date{\today}
\maketitle
\begin{abstract}
Coefficients in the expansions of the form
\begin{displaymath}
\frac{\partial P_{n}(\lambda;z)}{\partial\lambda}
=\sum_{k=0}^{n}a_{nk}(\lambda)P_{k}(\lambda;z),
\end{displaymath}
where $P_{n}(\lambda;z)$ is the $n$th classical (the generalized
Laguerre, Gegenbauer or Jacobi) orthogonal polynomial of variable
$z$ and $\lambda$ is a parameter, are evaluated. A method we
adopt in the present paper differs from that used by Fr{\"o}hlich
[Integral Transforms Spec.\ Funct.\ 2 (1994) 253] for the Jacobi
polynomials and by Koepf [Integral Transforms Spec.\ Funct.\ 5 (1997)
69] for the generalized Laguerre and the Gegenbauer polynomials.
\vskip1ex
\noindent
\textbf{KEY WORDS:} orthogonal polynomials; 
generalized Laguerre polynomials; Gegenbauer polynomials;
Jacobi polynomials; parameter derivatives
\vskip1ex
\noindent
\textbf{AMS subject classification:} 42C05, 33C45, 42C10
\end{abstract}
\maketitle
%
%
\section{Introduction}
\label{I}
\setcounter{equation}{0}
The problem of evaluation of coefficients in the expansions of the
form
\begin{equation}
\frac{\partial P_{n}(\lambda;z)}{\partial\lambda}
=\sum_{k=0}^{n}a_{nk}(\lambda)P_{k}(\lambda;z),
\label{1.1}
\end{equation}
where $P_{n}(\lambda;z)$ is the $n$th classical (the generalized
Laguerre, Gegenbauer or Jacobi \cite{Szeg39,Magn66,Hoch72})
orthogonal polynomial of variable $z$, with $\lambda$ being a
parameter, is not new. Fr{\"o}hlich \cite{Froh94} investigated the
case of the Jacobi polynomials, 
$P_{n}(\lambda;z)=P_{n}^{(\alpha,\beta)}(z)$, with either of the two
parameters $\alpha$ or $\beta$ playing the role of $\lambda$. He
started with a relationship between the Jacobi polynomials and the
Gauss hypergeometric function, found partial derivatives of the
latter with respect to the second and the third parameters, and then
ingeniously used some combinatorial identities to obtain closed-form
expressions for the coefficients in two (one corresponding to the
choice $\lambda=\alpha$, the other to the choice $\lambda=\beta$)
relevant expansions of the form (\ref{1.1}). Using a notation more
suitable for the present purposes than that used originally in Ref.\
\cite{Froh94}, Fr{\"o}hlich's results may be written as (cf Ref.\
\cite[Eqs.\ (1.2.23.1) and (1.2.23.2)]{Bryc08})
\begin{eqnarray}
\frac{\partial P_{n}^{(\alpha,\beta)}(z)}{\partial\alpha}
&=& [\psi(2n+\alpha+\beta+1)-\psi(n+\alpha+\beta+1)]
P_{n}^{(\alpha,\beta)}(z)
\nonumber \\
&& +\,\frac{\Gamma(n+\beta+1)}{\Gamma(n+\alpha+\beta+1)}
\sum_{k=0}^{n-1}\frac{2k+\alpha+\beta+1}{(n-k)(k+n+\alpha+\beta+1)}
\nonumber \\
&& \quad \times\frac{\Gamma(k+\alpha+\beta+1)}{\Gamma(k+\beta+1)}
P_{k}^{(\alpha,\beta)}(z)
\label{1.2}
\end{eqnarray}
and
\begin{eqnarray}
\frac{\partial P_{n}^{(\alpha,\beta)}(z)}{\partial\beta}
&=& [\psi(2n+\alpha+\beta+1)-\psi(n+\alpha+\beta+1)]
P_{n}^{(\alpha,\beta)}(z)
\nonumber \\
&& +\,(-)^{n}\frac{\Gamma(n+\alpha+1)}{\Gamma(n+\alpha+\beta+1)}
\sum_{k=0}^{n-1}(-)^{k}
\frac{2k+\alpha+\beta+1}{(n-k)(k+n+\alpha+\beta+1)}
\nonumber \\
&& \quad \times\frac{\Gamma(k+\alpha+\beta+1)}{\Gamma(k+\alpha+1)}
P_{k}^{(\alpha,\beta)}(z),
\label{1.3}
\end{eqnarray}
where
\begin{equation}
\psi(\zeta)=\frac{1}{\Gamma(\zeta)}
\frac{\mathrm{d}\Gamma(\zeta)}{\mathrm{d}\zeta}
\label{1.4}
\end{equation}
is the digamma function. Somewhat later, Koepf \cite{Koep97} pointed
out that one might use the Fr{\"o}hlich's results (\ref{1.2}) and
(\ref{1.3}) and the following known relationship between the
Gegenbauer and the Jacobi polynomials:
\begin{equation}
C_{n}^{(\lambda)}(z)
=\frac{\Gamma(\lambda+\frac{1}{2})}{\Gamma(2\lambda)}
\frac{\Gamma(n+2\lambda)}{\Gamma(n+\lambda+\frac{1}{2})}
P_{n}^{(\lambda-1/2,\lambda-1/2)}(z)
\label{1.5}
\end{equation}
to obtain an expansion of the form (\ref{1.1}) for
$C_{n}^{(\lambda)}(z)$. With some simplifying notational changes,
Koepf's result is (cf Ref.\ \cite[Eq.\ (1.22.2.2)]{Bryc08})
\begin{equation}
\frac{\partial C_{n}^{(\lambda)}(z)}{\partial\lambda}
=[\psi(n+\lambda)-\psi(\lambda)]C_{n}^{(\lambda)}(z)
+4\sum_{k=0}^{n-1}\frac{1+(-)^{k+n}}{2}\frac{k+\lambda}
{(n-k)(k+n+2\lambda)}C_{k}^{(\lambda)}(z),
\label{1.6}
\end{equation}
which may be also rewritten as (cf Ref.\ \cite[Eq.\
(1.22.2.1)]{Bryc08})
\begin{equation}
\frac{\partial C_{n}^{(\lambda)}(z)}{\partial\lambda}
=[\psi(n+\lambda)-\psi(\lambda)]C_{n}^{(\lambda)}(z)
+\sum_{k=1}^{\mathrm{int}(n/2)}\frac{n-2k+\lambda}
{k(n-k+\lambda)}C_{n-2k}^{(\lambda)}(z).
\label{1.7}
\end{equation}
Moreover, Koepf \cite{Koep97} observed that combining the
Fr{\"o}hlich's result (\ref{1.3}) with the known formula
\begin{equation}
L_{n}^{(\alpha)}(z)=\lim_{\beta\to\infty}P_{n}^{(\alpha,\beta)}
\left(1-\frac{2z}{\beta}\right),
\label{1.8}
\end{equation}
relating the generalized Laguerre and the Jacobi polynomials, one
obtains the remarkably simple expansion
\begin{equation}
\frac{\partial L_{n}^{(\lambda)}(z)}{\partial\lambda}
=\sum_{k=0}^{n-1}\frac{1}{n-k}L_{k}^{(\lambda)}(z),
\label{1.9}
\end{equation}
which finds applications in relativistic quantum mechanics, where it
is used to determine corrections of the order $\alpha^{2}$ ($\alpha$
is the fine-structure constant) to wave functions of particles bound
in the Coulomb potential. In Ref.\ \cite{Koep97}, Koepf gave an
alternative proof of this formula.

To make this brief survey complete, we have to mention that in Ref.\
\cite{Ronv00} Ronveaux \emph{et al.\/} found recurrence relations for
coefficients in an expansion more general than the one in Eq.\
(\ref{1.1}), namely
\begin{equation}
\frac{\partial^{m}P_{n}(\lambda;z)}{\partial\lambda^{m}}
=\sum_{k=0}^{n}a_{nk}(m,\lambda)P_{k}(\lambda;z)
\qquad (m\in\mathbb{N}),
\label{1.10}
\end{equation}
and in some particular cases they were able to obtain the
coefficients $a_{nk}(m,\lambda)$ in closed forms.

Our interest in the matter discussed in the present paper arose and
evolved in the course of carrying out research on derivatives of the
associated Legendre function of the first kind, $P_{\nu}^{\mu}(z)$,
with respect to its parameters $\mu$ and $\nu$
\cite{Szmy09a,Szmy09b,Szmy09c}. In particular, we have found
\cite{Szmy09b} that the derivatives
$[\partial P_{\nu}^{\pm M}(z)/\partial\nu]_{\nu=N}$, with
$M,N\in\mathbb{N}$, may be expressed in terms of the parameter
derivatives of the Jacobi polynomials, $\partial
P_{n}^{(\alpha,\beta)}(z)/\partial\beta$, with suitably chosen $n$,
$\alpha$ and $\beta$. Being at the time of doing the aforementioned
research unaware of the Fr{\"o}hlich's paper \cite{Froh94}, we
derived the relationships (\ref{1.2}) and (\ref{1.3}) independently.
Then it appeared to us that our way of reasoning might be applied to
derive the relationships (\ref{1.6}) and (\ref{1.9}) as well. Since
our approach turns out to be entirely different from the methods
adopted by Fr{\"o}hlich \cite{Froh94} and Koepf \cite{Koep97} (as
well as from the one used by Ronveaux \emph{et al.\/} \cite{Ronv00}),
we believe it is worthy to be presented.

Thus, in the present paper, we re-derive the expansions (\ref{1.2}),
(\ref{1.3}), (\ref{1.6}) and (\ref{1.9}). The path we shall follow
for each particular polynomial family comprises the following steps.
First, using the lemma proved in the Appendix and explicit
representations of the polynomials in terms of powers of
$z-z_{0}$, with suitably chosen $z_{0}$, for each family we shall
deduce the coefficient $a_{nn}(\lambda)$ standing in Eq.\ (\ref{1.1})
at $P_{n}(\lambda;z)$. Then, from the homogeneous differential
equation satisfied by $P_{n}(\lambda;z)$, in each particular case we
shall derive an inhomogeneous differential equation obeyed by
$\partial P_{n}(\lambda;z)/\partial\lambda$. Next, we shall show how,
in each case, after exploiting a relevant Christoffel--Darboux
identity, the inhomogeneity in that differential equation may be cast
to a form of a linear combination of the polynomials
$P_{k}(\lambda;z)$ with $0\leqslant k\leqslant n-1$. In the final
step, we shall insert this combination and the expansion (\ref{1.1})
into the aforementioned inhomogeneous differential equation for
$\partial P_{n}(\lambda;z)/\partial\lambda$, and then
straightforwardly deduce the coefficients $a_{nk}(\lambda)$ with
$0\leqslant k\leqslant n-1$.

All properties of the classical orthogonal polynomials we have
exploited in the present paper may be found in Ref.\ \cite{Hoch72}.
%
%
\section{The generalized Laguerre polynomials}
\label{II}
\setcounter{equation}{0}
The explicit formula defining the generalized Laguerre polynomials is
\begin{equation}
L_{n}^{(\lambda)}(z)=\sum_{k=0}^{n}\frac{(-)^{k}}{k!}
{n+\lambda \choose n-k} z^{k}.
\label{2.1}
\end{equation}
The coefficient at $z^{n}$ is
\begin{equation}
k_{n}(\lambda)=\frac{(-)^{n}}{n!}
\label{2.2}
\end{equation}
and is seen to be independent of $\lambda$. Hence, using Eq.\
(\ref{A.1}), we deduce that the coefficient $a_{nn}(\lambda)$ in the
expansion
\begin{equation}
\frac{\partial L_{n}^{(\lambda)}(z)}{\partial\lambda}
=\sum_{k=0}^{n}a_{nk}(\lambda)L_{k}^{(\lambda)}(z)
\label{2.3}
\end{equation}
vanishes:
\begin{equation}
a_{nn}(\lambda)=0.
\label{2.4}
\end{equation}

The coefficients $a_{nk}(\lambda)$ with $0\leqslant k\leqslant n-1$
may be found through the following reasoning. The differential
equation obeyed by $L_{n}^{(\lambda)}(z)$ is
\begin{equation}
\left[z\frac{\mathrm{d}^{2}}{\mathrm{d}z^{2}}
+(\lambda+1-z)\frac{\mathrm{d}}{\mathrm{d}z}+n\right]
L_{n}^{(\lambda)}(z)=0,
\label{2.5}
\end{equation}
which implies that
\begin{equation}
\left[z\frac{\mathrm{d}^{2}}{\mathrm{d}z^{2}}
+(\lambda+1-z)\frac{\mathrm{d}}{\mathrm{d}z}+n\right]
\frac{\partial L_{n}^{(\lambda)}(z)}{\partial\lambda}
=-\frac{\mathrm{d}L_{n}^{(\lambda)}(z)}{\mathrm{d}z}.
\label{2.6}
\end{equation}
We shall express the derivative
$\mathrm{d}L_{n}^{(\lambda)}(z)/\mathrm{d}z$ as a linear combination
of the polynomials $L_{k}^{(\lambda)}(z)$ with $0\leqslant k\leqslant
n-1$. To this end, at first we invoke the known relation
\begin{equation}
\frac{\mathrm{d}L_{n}^{(\lambda)}(z)}{\mathrm{d}z}
=\frac{nL_{n}^{(\lambda)}(z)-(n+\lambda)L_{n-1}^{(\lambda)}(z)}{z}.
\label{2.7}
\end{equation}
Next, we recall that the pertinent Christoffel--Darboux identity is
\begin{equation}
\sum_{k=0}^{n-1}\frac{k!}{\Gamma(k+\lambda+1)}
L_{k}^{(\lambda)}(z)L_{k}^{(\lambda)}(z')
=-\frac{n!}{\Gamma(n+\lambda)}
\frac{L_{n}^{(\lambda)}(z)L_{n-1}^{(\lambda)}(z')
-L_{n-1}^{(\lambda)}(z)L_{n}^{(\lambda)}(z')}{z-z'}.
\label{2.8}
\end{equation}
If we set here $z'=0$ and use
\begin{equation}
L_{k}^{(\lambda)}(0)=\frac{\Gamma(k+\lambda+1)}{k!\Gamma(\lambda+1)},
\label{2.9}
\end{equation}
Eq.\ (\ref{2.8}) becomes
\begin{equation}
\sum_{k=0}^{n-1}L_{k}^{(\lambda)}(z)
=-\frac{nL_{n}^{(\lambda)}(z)-(n+\lambda)L_{n-1}^{(\lambda)}(z)}{z}.
\label{2.10}
\end{equation}
Combining Eqs.\ (\ref{2.7}) and (\ref{2.10}) yields
\begin{equation}
\frac{\mathrm{d}L_{n}^{(\lambda)}(z)}{\mathrm{d}z}
=-\sum_{k=0}^{n-1}L_{k}^{(\lambda)}(z).
\label{2.11}
\end{equation}
Next, we insert Eqs.\ (\ref{2.3}) and (\ref{2.11}) into the left- and
right-hand sides of Eq.\ (\ref{2.6}), respectively, and then simplify
the left-hand side of the resulting equation with the aid of Eq.\
(\ref{2.5}). In this way, we obtain
\begin{equation}
\sum_{k=0}^{n-1}(n-k)a_{nk}(\lambda)L_{k}^{(\lambda)}(z)
=\sum_{k=0}^{n-1}L_{k}^{(\lambda)}(z).
\label{2.12}
\end{equation}
Equating coefficients at $L_{k}^{(\lambda)}(z)$ on both sides of Eq.\
(\ref{2.12}), we arrive at the sought expression for
$a_{nk}(\lambda)$:
\begin{equation}
a_{nk}(\lambda)=\frac{1}{n-k}
\qquad (0\leqslant k\leqslant n-1).
\label{2.13}
\end{equation}

From Eqs.\ (\ref{2.3}), (\ref{2.4}) and (\ref{2.13}) the relationship
(\ref{1.9}) follows.
%
%
\section{The Gegenbauer polynomials}
\label{III}
\setcounter{equation}{0}
The Gegenbauer polynomials may be defined through the formula
\begin{equation}
C_{n}^{(\lambda)}(z)
=\frac{\Gamma(\lambda+\frac{1}{2})}{\Gamma(2\lambda)}
\sum_{k=0}^{n}\frac{\Gamma(k+n+2\lambda)}
{k!(n-k)!\Gamma(k+\lambda+\frac{1}{2})}
\left(\frac{z-1}{2}\right)^{k}.
\label{3.1}
\end{equation}

To find the coefficient $a_{nn}^{(\lambda)}$ in the expansion
\begin{equation}
\frac{\partial C_{n}^{(\lambda)}(z)}{\partial\lambda}
=\sum_{k=0}^{n}a_{nk}(\lambda)C_{k}^{(\lambda)}(z),
\label{3.2}
\end{equation}
we infer from Eq.\ (\ref{3.1}) that the coefficient at $(z-1)^{n}$ is
\begin{equation}
k_{n}(\lambda)=\frac{\Gamma(\lambda+\frac{1}{2})}{\Gamma(2\lambda)}
\frac{\Gamma(2n+2\lambda)}{2^{n}n!\Gamma(n+\lambda+\frac{1}{2})}
=\frac{2^{n}\Gamma(n+\lambda)}{n!\Gamma(\lambda)}.
\label{3.3}
\end{equation}
Hence, application of the lemma from the Appendix gives
\begin{equation}
a_{nn}(\lambda)=\psi(n+\lambda)-\psi(\lambda).
\label{3.4}
\end{equation}

To proceed further, we recall that the Gegenbauer polynomials obey
the differential identity
\begin{equation}
\left[(1-z^{2})\frac{\mathrm{d}^{2}}{\mathrm{d}z^{2}}
-(2\lambda+1)z\frac{\mathrm{d}}{\mathrm{d}z}+n(n+2\lambda)\right]
C_{n}^{(\lambda)}(z)=0.
\label{3.5}
\end{equation}
Differentiating with respect to $\lambda$, we obtain
\begin{equation}
\left[(1-z^{2})\frac{\mathrm{d}^{2}}{\mathrm{d}z^{2}}
-(2\lambda+1)z\frac{\mathrm{d}}{\mathrm{d}z}+n(n+2\lambda)\right]
\frac{\partial C_{n}^{(\lambda)}(z)}{\partial\lambda}
=2z\frac{\mathrm{d}C_{n}^{(\lambda)}(z)}{\mathrm{d}z}
-2nC_{n}^{(\lambda)}(z).
\label{3.6}
\end{equation}
The next step is a bit tricky. It consists of writing
\begin{equation}
2z\frac{\mathrm{d}C_{n}^{(\lambda)}(z)}{\mathrm{d}z}
-2nC_{n}^{(\lambda)}(z)
=\left[(z+1)\frac{\mathrm{d}C_{n}^{(\lambda)}(z)}{\mathrm{d}z}
-nC_{n}^{(\lambda)}(z)\right]
+\left[(z-1)\frac{\mathrm{d}C_{n}^{(\lambda)}(z)}{\mathrm{d}z}
-nC_{n}^{(\lambda)}(z)\right].
\label{3.7}
\end{equation}
Using the known relation
\begin{equation}
(z^{2}-1)\frac{\mathrm{d}C_{n}^{(\lambda)}(z)}{\mathrm{d}z}
=nzC_{n}^{(\lambda)}(z)-(n+2\lambda-1)C_{n-1}^{(\lambda)}(z),
\label{3.8}
\end{equation}
Eq.\ (\ref{3.7}) may be cast to the form
\begin{equation}
2z\frac{\mathrm{d}C_{n}^{(\lambda)}(z)}{\mathrm{d}z}
-2nC_{n}^{(\lambda)}(z)
=\frac{nC_{n}^{(\lambda)}(z)-(n+2\lambda-1)C_{n-1}^{(\lambda)}(z)}
{z-1}
-\frac{nC_{n}^{(\lambda)}(z)+(n+2\lambda-1)C_{n-1}^{(\lambda)}(z)}
{z+1}.
\label{3.9}
\end{equation}
Now, the Christoffel--Darboux identity for the Gegenbauer polynomials
is
\begin{equation}
\sum_{k=0}^{n-1}\frac{k!(k+\lambda)}{\Gamma(k+2\lambda)}
C_{k}^{(\lambda)}(z)C_{k}^{(\lambda)}(z')
=\frac{n!}{2\Gamma(n+2\lambda-1)}
\frac{C_{n}^{(\lambda)}(z)C_{n-1}^{(\lambda)}(z')
-C_{n-1}^{(\lambda)}(z)C_{n}^{(\lambda)}(z')}{z-z'}.
\label{3.10}
\end{equation}
If we set in this identity $z'=\pm1$ and use
\begin{equation}
C_{k}^{(\lambda)}(\pm1)=(\pm)^{k}
\frac{\Gamma(k+2\lambda)}{k!\Gamma(2\lambda)},
\label{3.11}
\end{equation}
we obtain
\begin{equation}
\sum_{k=0}^{n-1}(\pm)^{k}(k+\lambda)C_{k}^{(\lambda)}(z)
=\frac{(\pm)^{n}}{2}
\frac{\pm nC_{n}^{(\lambda)}(z)-(n+2\lambda-1)C_{n-1}^{(\lambda)}(z)}
{z\mp1}.
\label{3.12}
\end{equation}
Combining this result with Eq.\ (\ref{3.9}) yields the expansion
\begin{equation}
2z\frac{\mathrm{d}C_{n}^{(\lambda)}(z)}{\mathrm{d}z}
-2nC_{n}^{(\lambda)}(z)
=4\sum_{k=0}^{n-1}\frac{1+(-)^{k+n}}{2}
(k+\lambda)C_{k}^{(\lambda)}(z).
\label{3.13}
\end{equation}
If we plug the expansions (\ref{3.2}) and (\ref{3.13}) into Eq.\
(\ref{3.6}), use Eq.\ (\ref{3.5}) and equate coefficients at
$C_{k}^{(\lambda)}(z)$ on both sides of the resulting equation, this
gives
\begin{equation}
a_{nk}(\lambda)=4\frac{1+(-)^{k+n}}{2}
\frac{k+\lambda}{(n-k)(k+n+2\lambda)}
\qquad (0\leqslant k\leqslant n-1).
\label{3.14}
\end{equation}

Insertion of Eqs.\ (\ref{3.4}) and (\ref{3.14}) into Eq.\ (\ref{3.2})
results in the expansion (\ref{1.6}).
%
%
\section{The Jacobi polynomials}
\label{IV}
\setcounter{equation}{0}
Finally, we shall determine coefficients in the expansion
\begin{equation}
\frac{\partial P_{n}^{(\alpha,\beta)}(z)}{\partial\alpha}
=\sum_{k=0}^{n}a_{nk}^{(\beta)}(\alpha)P_{k}^{(\alpha,\beta)}(z),
\label{4.1}
\end{equation}
where 
\begin{equation}
P_{n}^{(\alpha,\beta)}(z)=\frac{\Gamma(n+\alpha+1)}
{\Gamma(n+\alpha+\beta+1)}
\sum_{k=0}^{n}\frac{\Gamma(k+n+\alpha+\beta+1)}
{k!(n-k)!\Gamma(k+\alpha+1)}\left(\frac{z-1}{2}\right)^{k}
\label{4.2}
\end{equation}
is the Jacobi polynomial. Observe that once these coefficients are
determined, coefficients in the counterpart expansion of $\partial
P_{n}^{(\alpha,\beta)}(z)/\partial\beta$ are also known, because in
virtue of the relationship
\begin{equation}
P_{n}^{(\alpha,\beta)}(z)=(-)^{n}P_{n}^{(\beta,\alpha)}(-z)
\label{4.3}
\end{equation}
one has
\begin{equation}
\frac{\partial P_{n}^{(\alpha,\beta)}(z)}{\partial\beta}
=\sum_{k=0}^{n}(-)^{k+n}a_{nk}^{(\alpha)}(\beta)
P_{k}^{(\alpha,\beta)}(z).
\label{4.4}
\end{equation}

The coefficient at $(z-1)^{n}$ on the right-hand side of Eq.\
(\ref{4.2}) is
\begin{equation}
k_{n}(\alpha,\beta)=\frac{\Gamma(2n+\alpha+\beta+1)}
{2^{n}n!\Gamma(n+\alpha+\beta+1)},
\label{4.5}
\end{equation}
hence, the coefficient $a_{nn}^{(\beta)}(\alpha)$ in the expansion
(\ref{4.1}) is
\begin{equation}
a_{nn}^{(\beta)}(\alpha)=\psi(2n+\alpha+\beta+1)
-\psi(n+\alpha+\beta+1).
\label{4.6}
\end{equation}

To find the remaining coefficients $a_{nk}^{(\beta)}(\alpha)$ with
$0\leqslant k\leqslant n-1$, we differentiate the identity
\begin{equation}
\left\{(1-z^{2})\frac{\mathrm{d}^{2}}{\mathrm{d}z^{2}}
+[\beta-\alpha-(\alpha+\beta+2)z]\frac{\mathrm{d}}{\mathrm{d}z}
+n(n+\alpha+\beta+1)\right\}P_{n}^{(\alpha,\beta)}(z)=0
\label{4.7}
\end{equation}
with respect to $\alpha$, obtaining
\begin{eqnarray}
&& \left\{(1-z^{2})\frac{\mathrm{d}^{2}}{\mathrm{d}z^{2}}
+[\beta-\alpha-(\alpha+\beta+2)z]\frac{\mathrm{d}}{\mathrm{d}z}
+n(n+\alpha+\beta+1)\right\}
\frac{\partial P_{n}^{(\alpha,\beta)}(z)}{\partial\alpha}
\nonumber \\
&& \hspace*{10em} 
=\,(z+1)\frac{\mathrm{d}P_{n}^{(\alpha,\beta)}(z)}{\mathrm{d}z}
-nP_{n}^{(\alpha,\beta)}(z).
\label{4.8}
\end{eqnarray}
Using the known relationship
\begin{eqnarray}
(2n+\alpha+\beta)(z^{2}-1)
\frac{\mathrm{d}P_{n}^{(\alpha,\beta)}(z)}{\mathrm{d}z}
&=& n[\beta-\alpha+(2n+\alpha+\beta)z]P_{n}^{(\alpha,\beta)}(z)
\nonumber \\
&& -\,2(n+\alpha)(n+\beta)P_{n-1}^{(\alpha,\beta)}(z),
\label{4.9}
\end{eqnarray}
the right-hand side of Eq.\ (\ref{4.8}) may be rewritten as
\begin{eqnarray}
&& (z+1)\frac{\mathrm{d}P_{n}^{(\alpha,\beta)}(z)}{\mathrm{d}z}
-nP_{n}^{(\alpha,\beta)}(z)=\frac{2(n+\beta)}{2n+\alpha+\beta}
\frac{nP_{n}^{(\alpha,\beta)}(z)
-(n+\alpha)P_{n-1}^{(\alpha,\beta)}(z)}{z-1}.
\label{4.10}
\end{eqnarray}
The Christoffel--Darboux formula for the Jacobi polynomials is
\begin{eqnarray}
&& \sum_{k=0}^{n-1}(2k+\alpha+\beta+1)
\frac{k!\Gamma(k+\alpha+\beta+1)}{\Gamma(k+\alpha+1)\Gamma(k+\beta+1)}
P_{k}^{(\alpha,\beta)}(z)P_{k}^{(\alpha,\beta)}(z')
\nonumber \\
&& \qquad\qquad
=\frac{2n!\Gamma(n+\alpha+\beta+1)}
{(2n+\alpha+\beta)\Gamma(n+\alpha)\Gamma(n+\beta)}
\frac{P_{n}^{(\alpha,\beta)}(z)P_{n-1}^{(\alpha,\beta)}(z')
-P_{n-1}^{(\alpha,\beta)}(z)P_{n}^{(\alpha,\beta)}(z')}{z-z'}.
\nonumber \\
&&
\label{4.11}
\end{eqnarray}
We set in this formula $z'=1$ and then use
\begin{equation}
P_{k}^{(\alpha,\beta)}(1)
=\frac{\Gamma(k+\alpha+1)}{k!\Gamma(\alpha+1)}.
\label{4.12}
\end{equation}
This yields
\begin{eqnarray}
&& \sum_{k=0}^{n-1}(2k+\alpha+\beta+1)
\frac{\Gamma(k+\alpha+\beta+1)}{\Gamma(k+\beta+1)}
P_{k}^{(\alpha,\beta)}(z)
\nonumber \\
&& \qquad\qquad =\frac{2\Gamma(n+\alpha+\beta+1)}
{(2n+\alpha+\beta)\Gamma(n+\beta)}
\frac{nP_{n}^{(\alpha,\beta)}(z)
-(n+\alpha)P_{n-1}^{(\alpha,\beta)}(z)}{z-1}.
\label{4.13}
\end{eqnarray}
Combining this result with Eq.\ (\ref{4.10}) gives us the following
expansion of the expression standing on the right-hand side of Eq.\
(\ref{4.8}):
\begin{eqnarray}
&& (z+1)\frac{\mathrm{d}P_{n}^{(\alpha,\beta)}(z)}{\mathrm{d}z}
-nP_{n}^{(\alpha,\beta)}(z)
\nonumber \\
&& \qquad\qquad
=\frac{\Gamma(n+\beta+1)}{\Gamma(n+\alpha+\beta+1)}
\sum_{k=0}^{n-1}(2k+\alpha+\beta+1)
\frac{\Gamma(k+\alpha+\beta+1)}{\Gamma(k+\beta+1)}
P_{k}^{(\alpha,\beta)}(z).
\label{4.14}
\end{eqnarray}
Substitution of Eq.\ (\ref{4.1}) into the left-hand side and of Eq.\
(\ref{4.14}) into the right-hand side of Eq.\ (\ref{4.8}), followed
by the use of Eq.\ (\ref{4.7}), after some obvious movements leads to
\begin{eqnarray}
a_{nk}^{(\beta)}(\alpha)
&=& \frac{2k+\alpha+\beta+1}{(n-k)(k+n+\alpha+\beta+1)}
\frac{\Gamma(n+\beta+1)\Gamma(k+\alpha+\beta+1)}
{\Gamma(k+\beta+1)\Gamma(n+\alpha+\beta+1)}
\nonumber \\
&& (0\leqslant k\leqslant n-1).
\label{4.15}
\end{eqnarray}

Plugging Eqs.\ (\ref{4.6}) and (\ref{4.15}) into Eq.\ (\ref{4.1})
results in the relation (\ref{1.2}), while Eqs.\ (\ref{4.4}),
(\ref{4.6}) and (\ref{4.15}) imply Eq.\ (\ref{1.3}).
%
%
\appendix
\section{Appendix}
\label{A}
\setcounter{equation}{0}
In this Appendix we shall prove the following
\begin{lem}
The coefficient $a_{nn}(\lambda)$ in the expansion (\ref{1.1})
is given by
\begin{equation}
a_{nn}(\lambda)=\frac{\partial\ln k_{n}(\lambda)}{\partial\lambda}, 
\label{A.1} 
\end{equation}
where $k_{n}(\lambda)$ is the coefficient at $(z-z_{0})^{n}$ (with
$z_{0}\in\mathbb{C}$ arbitrary) in the expansion of the polynomial
$P_{n}(\lambda;z)$ in powers of $z-z_{0}$.
\end{lem}
\begin{proof}
We have
\begin{equation}
P_{n}(\lambda;z)=k_{n}(\lambda)(z-z_{0})^{n}
+Q_{n-1}(\lambda,z_{0};z)
\label{A.2}
\end{equation}
(obviously, $k_{n}(\lambda)$ is independent of $z_{0}$), where 
$Q_{n-1}(\lambda,z_{0};z)$ is a polynomial in $z-z_{0}$ of degree
$n-1$. Hence, it follows that
\begin{equation}
\frac{\partial P_{n}(\lambda;z)}{\partial\lambda}
=\frac{\partial k_{n}(\lambda)}{\partial\lambda}(z-z_{0})^{n}
+\frac{\partial Q_{n-1}(\lambda,z_{0};z)}{\partial\lambda}.
\label{A.3}
\end{equation}
On the other side, from Eqs.\ (\ref{1.1}) and (\ref{A.2}) we have
\begin{eqnarray}
\frac{\partial P_{n}(\lambda;z)}{\partial\lambda}
&=& a_{nn}(\lambda)P_{n}(\lambda;z)
+\sum_{k=0}^{n-1}a_{nk}(\lambda)P_{k}(\lambda;z)
\nonumber \\
&=& a_{nn}(\lambda)k_{n}(\lambda)(z-z_{0})^{n}
+R_{n-1}(\lambda,z_{0};z),
\label{A.4}
\end{eqnarray}
where $R_{n-1}(\lambda,z_{0};z)$ is a polynomial in $z-z_{0}$ of
degree $n-1$. Equating coefficients at $(z-z_{0})^{n}$ on the
right-hand sides of Eqs.\ (\ref{A.3}) and (\ref{A.4}) yields the
expression (\ref{A.1}).
\end{proof}
%
%

%
\end{document}